\documentclass{article}
\usepackage{kz}
\usepackage{lmodern}
\usepackage{subcaption}

\DeclareMathOperator{\rank}{rank}
\DeclareMathOperator{\Supp}{Supp}
\DeclareMathOperator{\TV}{TV}

\title{Binary Expansion Group Intersection Network}
\author{Sicheng Zhou\footnote{Sicheng Zhou is an undergraduate student (E-mail: sichengz@mit.edu), Department of Electrical Engineering and Computer Science. Massachusetts Institute of Technology, Cambridge, MA 02139.} and Kai Zhang\footnote{Kai Zhang is a Professor (Corresponding author. E-mail: zhangk@email.unc.edu), Department of Statistics and Operations Research, University of North Carolina at Chapel Hill, Chapel Hill, NC 27599.}}
\begin{document}
\maketitle\large

\begin{abstract}
Conditional independence is central to modern statistics, but beyond special parametric families it rarely admits an exact covariance characterization. We introduce the \emph{binary expansion group intersection network} (BEGIN), a distribution-free graphical representation of conditional independence for multivariate binary data and bit-encoded multinomial variables. For arbitrary binary random vectors and bit representations of multinomial variables, we prove that conditional independence is equivalent to a sparse linear representation of conditional expectations, to a block factorization of the corresponding interaction covariance matrix, and to block diagonality of an associated generalized Schur complement. The resulting graph is indexed by the intersection of multiplicative groups of binary interactions, yielding an analogue of Gaussian graphical modeling beyond the Gaussian setting. This viewpoint treats data bits as atoms and local BEGIN molecules as building blocks for large Markov random fields. We also show how dyadic bit representations allow BEGIN to approximate conditional independence for general random vectors under mild regularity conditions. A key technical device is the \emph{Hadamard prism}, a linear map that links interaction covariances to group structure.  
\end{abstract}

\section{Introduction}
Conditional independence is a cornerstone of statistical reasoning. It underlies the interpretation of multivariate associations, the construction of graphical models, and many procedures in causal inference, variable selection, and structure learning. In classical settings, conditional independence is often studied through parametric models whose covariance structure has a direct probabilistic interpretation. In many modern applications, however, the available data are heterogeneous, high dimensional, or only weakly modeled, so exact parametric assumptions are difficult to justify.

This tension is especially acute in distribution-free inference. On the one hand, conditional independence is one of the most natural ways to formalize ``no direct association after adjustment.'' On the other hand, fully distribution-free conditional-independence testing is fundamentally impossible without additional structure; see, for example, \citet{shah2020hardness}. The challenge, then, is to identify structure that is both mathematically exact and as assumption-lean as possible; see, for example, \citet{christgau2023nonparametric}.

This paper builds on the multiresolution viewpoint of \emph{binary expansion statistics} \citep{zhang2019bet, zhang2021beauty, brown2025belief}, which treats data bits as atomic units of information. At the bit level, binary variables exhibit an exact linearity that does not persist for general variables. That observation suggests asking whether conditional independence can be characterized exactly through covariances once one works with suitable binary interaction features beyond the original variables themselves.

Graphical models provide the natural language for this question \citep{lauritzen1996graphical, wainwright2008graphical, koller2009probabilistic, drton2017structure}. In Gaussian graphical models, conditional independence is equivalent to sparsity of the precision matrix. That equivalence drives both the interpretation of the graph and the design of scalable estimation procedures. Outside the Gaussian family, however, zeros in the inverse covariance matrix generally do \emph{not} encode conditional independence.

For discrete data, several important alternatives are available. Classical log-linear models already provide exact factorization-based characterizations of binary conditional independence, whereas Ising models and more general Markov random fields impose explicit factorization assumptions. \citet{loh2012structure} study generalized covariance matrices built from sufficient statistics of multinomial exponential-family models and show that their inverses can reflect graph structure under those modeling assumptions. \citet{lauritzen2021total} derive strong implications for binary distributions under total positivity constraints. Our goal here is different: an exact covariance-based characterization in an interaction basis for arbitrary multivariate binary distributions.

With the central result Theorem~\ref{thm:begin}, this paper makes the following four main contributions:
\begin{enumerate}[(i)]
\item We prove that for arbitrary binary random vectors $(\bA,\bB,\bC)$, conditional independence $\bA \independent \bC \mid \bB$ is equivalent to a covariance structure indexed by the intersection of multiplicative groups of binary interactions. The key object is not the inverse of the full covariance matrix, but a generalized Schur complement associated with the interaction block generated by $\bB$.
\item We show that the same framework remains valid for bit-encoded multinomial variables, including rank-deficient cases created by deterministic constraints or structural zeros.
\item We introduce the \emph{Hadamard prism}, a convenient linear map for the covariance algebra of binary interactions that clarifies the link between covariance of binary interactions, Walsh--Hadamard transforms, and Boolean Fourier analysis.
\item We extend the framework beyond discrete data by showing that dyadic quantizations preserve conditional independence asymptotically and yield explicit approximation bounds under H\"older-type continuity of the relevant conditional laws.
\end{enumerate}

Taken together, these results yield the graph interpretation in Corollary~\ref{cor:begin}, which continues to hold in singular multinomial encodings. To our knowledge, this corollary provides the first exact distribution-free covariance-based graphical characterization of conditional independence for arbitrary multivariate binary data. The graph that emerges is indexed not merely by the original variables but by intersections of multiplicative groups generated by their interactions. For that reason, we call the resulting representation the \emph{binary expansion group intersection network} (BEGIN). 

A useful way to position BEGIN is relative to Gaussian graphical models and generalized covariance constructions. BEGIN is Gaussian-like in spirit because conditional independence is read off from a sparse matrix object and the resulting structure suggests nodewise projection viewpoints. It is fundamentally non-Gaussian, however, because the relevant nodes are interaction features and the correct matrix object is a generalized Schur complement rather than an ordinary precision matrix. Compared with exponential-family generalized covariance methods, BEGIN requires weaker modeling assumptions: it does not rely on strict positivity, a prescribed clique factorization, or a particular parametric family. By viewing data bits as atoms and BEGIN as molecules, we provide examples showing how local BEGIN structures can serve as building blocks for larger Markov random fields

Section~\ref{sec:begin} develops the BEGIN characterization for binary variables and bit-represented multinomial variables. Section~\ref{sec:approx} studies dyadic approximations for general random vectors. Section~\ref{sec:discussion} closes with brief remarks on implications and future work. Proofs are provided in the supplementary material.

\section*{Notation}
We use the following notation throughout. For $\by \in \RR^p$, $\diag(\by)$ denotes the diagonal matrix with diagonal entries $\by$. For a matrix $\Mb$, $\Mb^{+}$ denotes the Moore--Penrose inverse. The symbol $\Hb_p$ denotes the $2^p \times 2^p$ Hadamard matrix obtained by Sylvester's construction,
\begin{equation}
\Hb_p=
\underbrace{
\begin{pmatrix}
1 & 1\\
1 & -1
\end{pmatrix}
\otimes \cdots \otimes
\begin{pmatrix}
1 & 1\\
1 & -1
\end{pmatrix}}_{p\text{ factors}}.
\end{equation}
For a binary random vector $\bX=(X_1,\ldots,X_p) \in \{\pm 1\}^p$, let $\langle \bX \rangle$ denote the multiplicative group generated by the coordinates of $\bX$, and for multiple binary vectors $\bX_1,\ldots,\bX_k$, let $\langle \bX_1,\ldots,\bX_k \rangle$ denote the group generated by the union of their coordinates. The \emph{power vector} $\bX^{\otimes} \in \{\pm 1\}^{2^p}$ collects the elements of $\langle \bX \rangle$ via
\begin{equation}\label{eq:power}
\bX^{\otimes}:=
\begin{pmatrix}1\\ X_1\end{pmatrix}
\otimes
\begin{pmatrix}1\\ X_2\end{pmatrix}
\otimes \cdots \otimes
\begin{pmatrix}1\\ X_p\end{pmatrix}.
\end{equation}
We index the coordinates of $\bX^{\otimes}$ by $\Lambda \in \{0,1\}^p$ and write $X_{\Lambda}:=\prod_{j=1}^p X_j^{\Lambda_j}$. For a matrix $\Mb$ and index sets $\cR,\cC$, $\Mb[\cR,\cC]$ denotes the corresponding submatrix. For a random object $\bZ$, $\Supp(\bZ)$ denotes its support with positive probabilities. For a set $\cS$, $|\cS|$ denotes its cardinality. For probability measures $P$ and $Q$ on a common measurable space, $\TV(P,Q):=\sup_{S}|P(S)-Q(S)|$ denotes total variation distance.

\section{Main Results}\label{sec:begin}
\subsection{Characterization of Conditional Independence for Binary Variables}
Let $\bA \in \{\pm 1\}^r$, $\bB \in \{\pm 1\}^s$, and $\bC \in \{\pm 1\}^t$ be binary random vectors.

\begin{definition}
We say that $\bA$ and $\bC$ are \emph{conditionally independent given $\bB$}, written $\bA \independent \bC \mid \bB$, if for all $\ba,\bb,\bc$ with $\Pb(\bB=\bb)>0$,
\[
\Pb(\bA=\ba,\bC=\bc \mid \bB=\bb)
=
\Pb(\bA=\ba \mid \bB=\bb)\,
\Pb(\bC=\bc \mid \bB=\bb).
\]
\end{definition}

To connect conditional independence with covariance structure, we allow zero-probability cells in the joint distribution of $(\bA,\bB,\bC)$. This is essential for multinomial variables encoded by data bits. If a multinomial variable has $m$ categories with positive probability and $2^{k-1}<m\le 2^k$, then it can be represented by $k$ bits, but that representation may impose deterministic constraints among interaction features and hence lead to singular covariance matrices. The following rank characterization makes that singularity transparent.

\begin{theorem}\label{thm:rank}
For a binary random vector $\bX \in \{\pm 1\}^p$, let
\[
\bSigma_{\langle \bX \rangle \setminus \{1\}}
:=
\Cov\bigl(\langle \bX \rangle \setminus \{1\},\, \langle \bX \rangle \setminus \{1\}\bigr)
\]
be the covariance matrix of the nonconstant elements of $\langle \bX \rangle$. Then
\[
\rank\bigl(\bSigma_{\langle \bX \rangle \setminus \{1\}}\bigr)=|\Supp(\bX)|-1.
\]
\end{theorem}

We now define the interaction index sets that govern BEGIN:
\[
\cB:=\langle \bB \rangle \setminus \{1\},
\qquad
\cL:=\langle \bA,\bB \rangle \setminus \langle \bB \rangle,
\qquad
\cR:=\langle \bB,\bC \rangle \setminus \langle \bB \rangle.
\]
Let $\bSigma$ be the covariance matrix of the concatenated interaction vector indexed by $\cB \cup \cL \cup \cR$, ordered as $(\cB,\cL,\cR)$. When the joint pmf of $(\bA,\bB,\bC)$ is strictly positive, $\bSigma$ is positive definite; in general it is only positive semidefinite.

\begin{theorem}\label{thm:begin}\emph{(BEGIN)}
The following statements are equivalent.
\begin{enumerate}[(a)]
\item $\bA \independent \bC \mid \bB$.

\item \emph{Sparse conditional-expectation representation.} For every $\Lambda_1 \in \{0,1\}^r$ and $\Lambda_2 \in \{0,1\}^t$, there exist coefficient vectors $\balpha_{\Lambda_1},\bgamma_{\Lambda_2} \in \RR^{2^s}$ such that
\begin{equation}
\begin{split}
\Eb\bigl[A_{\Lambda_1} \mid \bB,\bC\bigr]
&=
\Eb\bigl[A_{\Lambda_1} \mid \bB\bigr]
=
\balpha_{\Lambda_1}^{\top}\bB^{\otimes},\\
\Eb\bigl[C_{\Lambda_2} \mid \bA,\bB\bigr]
&=
\Eb\bigl[C_{\Lambda_2} \mid \bB\bigr]
=
\bgamma_{\Lambda_2}^{\top}\bB^{\otimes}.
\end{split}
\end{equation}

\item \emph{Block factorization of covariance blocks.} There exist matrices $\Mb_1$ and $\Mb_2$ such that
\begin{equation}
\bSigma=
\begin{pmatrix}
\bSigma_{\cB} & \bSigma_{\cB}\Mb_1^{\top} & \bSigma_{\cB}\Mb_2^{\top}\\
\Mb_1\bSigma_{\cB} & \bSigma_{\cL} & \Mb_1\bSigma_{\cB}\Mb_2^{\top}\\
\Mb_2\bSigma_{\cB} & \Mb_2\bSigma_{\cB}\Mb_1^{\top} & \bSigma_{\cR}
\end{pmatrix}.
\end{equation}

\item \emph{Block-diagonal generalized Schur complement.} The generalized Schur complement of $\bSigma_{\cB}$ in $\bSigma$,
\[
\Sbb
:=
\bSigma[\cL \cup \cR,\, \cL \cup \cR]
-
\bSigma[\cL \cup \cR,\, \cB]\,
\bSigma_{\cB}^{+}\,
\bSigma[\cB,\, \cL \cup \cR],
\]
is block diagonal with respect to $(\cL,\cR)$; that is,
\[
\Sbb=
\begin{pmatrix}
\Sbb_{\cL} & \zero\\
\zero & \Sbb_{\cR}
\end{pmatrix}.
\]
\end{enumerate}
\end{theorem}

Theorem~\ref{thm:begin} identifies the exact covariance object behind conditional independence for binary data. Part~(b) follows from the binary expansion linear effect (BELIEF) representation in \citet{brown2025belief}. The novelty here is that this interaction-level linearity is equivalent to the covariance factorization in part~(c) and to the generalized Schur-complement sparsity in part~(d). Under conditional independence, the conditional expectation of every interaction on the left or right depends on $(\bA,\bB,\bC)$ only through the $2^s$ interaction coordinates in $\bB^{\otimes}$; without conditional independence, the corresponding representations generally require $2^{r+s}$ or $2^{s+t}$ coefficients.

Parts (c) and (d) of Theorem~\ref{thm:begin} provide a one-to-one characterization of conditional independence in terms of the covariance structure of $\langle \bA,\bB \rangle \cap \langle \bB,\bC \rangle$ for an arbitrary binary vector $(\bA,\bB,\bC)$. We emphasize that conditional independence in binary variables must be expressed through intersections of groups: for binary variables, the $\sigma$-field is determined by the group they generate. If one replaces $\cB$ by a proper subset of the group intersection, the equivalence can fail in either direction: sparsity of the corresponding Schur complement need not imply conditional independence, and conditional independence need not imply sparsity. The supplementary material provides explicit counterexamples for both failures. Because the relevant covariance blocks are indexed by intersections of multiplicative groups of binary interactions, we call this structure the \emph{binary expansion group intersection network} (BEGIN).

Theorem~\ref{thm:begin} also shows why the ordinary inverse covariance matrix is not the right object for describing conditional independence outside the Gaussian setting. When $\bSigma$ is singular, the Moore--Penrose inverse $\bSigma^{+}$ need not reflect the relevant sparsity pattern. BEGIN instead isolates the interaction block generated by $\bB$ and the associated generalized Schur complement. This is motivated by Theorem~2.5 of \citet{brown2025belief}, which implies that the rows of $\bSigma[\cL\cup\cR,\cB]$ lie in the row space of $\bSigma_{\cB}$. Equivalently, if we define
\begin{equation}\label{eq:rowspace_general}
\Mb:=\bSigma[\cL\cup\cR,\cB] \bSigma_{\cB}^{+},
\end{equation}
then
\begin{equation}\label{eq:rowspace_factor}
\bSigma[\cL\cup\cR,\cB]=\Mb\bSigma_{\cB},
\qquad
\bSigma[\cB,\cL\cup\cR]=\bSigma_{\cB}\Mb^{\top}.
\end{equation}
This deterministic row-space identity motivates the use of the Schur--Banachiewicz inverse of \cite{ouellette1981schur}, rather than the more commonly used Moore--Penrose inverse $\bSigma^{+}$.

\begin{definition}
Define the \emph{Schur--Banachiewicz generalized inverse} of $\bSigma$ by
\begin{equation}\label{eq:sbinv}
\bOmega:=
\begin{pmatrix}
\bSigma_{\cB}^{+}+\bSigma_{\cB}^{+}\Fb\Sbb^{+}\Fb^{\top}\bSigma_{\cB}^{+} & -\bSigma_{\cB}^{+}\Fb\Sbb^{+}\\[2pt]
-\Sbb^{+}\Fb^{\top}\bSigma_{\cB}^{+} & \Sbb^{+}
\end{pmatrix},
\qquad
\Fb:=\bSigma[\cB,\cL\cup\cR].
\end{equation}
\end{definition}

We use the Schur--Banachiewicz inverse rather than $\bSigma^{+}$ because its $\bOmega[\cL\cup\cR,\cL\cup\cR]$ is exactly $\Sbb^{+}$. It therefore preserves the separation structure in Corollary~\ref{cor:begin} induced by the generalized Schur complement, even when $\bSigma$ is singular, and yields an exact characterization for rank-deficient multinomial encodings.

\begin{corollary}\label{cor:begin}
The Schur--Banachiewicz inverse $\bOmega$ is symmetric and satisfies $\bSigma\,\bOmega\,\bSigma=\bSigma.$
Moreover, the following statements are equivalent.
\begin{enumerate}[(a)]
\item $\bA \independent \bC \mid \bB$.
\item $\bOmega[\cL,\cR]=\zero$.
\item In the undirected graph with vertex set $\cB\cup\cL\cup\cR$ and an edge between two distinct vertices whenever the corresponding entry of $\bOmega$ is nonzero, the set $\cB$ separates $\cL$ from $\cR$.
\end{enumerate}
\end{corollary}

Corollary~\ref{cor:begin} shows that BEGIN plays the same structural role for binary interaction features that the precision matrix plays in Gaussian graphical models. In particular, since $\bOmega[\cL\cup\cR,\cL\cup\cR]=\Sbb^{+}$, the graph can be read directly from the sparsity pattern of the Schur--Banachiewicz inverse. As in Gaussian graphical modeling, this viewpoint suggests possible nodewise estimation strategies, though developing their finite-sample theory is beyond the scope of this note. Unlike the Gaussian case, however, the relevant nodes are interaction features and the underlying matrix may be singular.

It is also useful to compare BEGIN with the generalized covariance approach of \citet{loh2012structure} and with the Ising model. A pairwise Ising model is a strictly positive Markov random field on the original variables, so its graph is specified by a factorized likelihood; classical log-linear models likewise provide exact factorization-based characterizations of binary conditional independence. \citet{loh2012structure} remain within discrete exponential-family graphical models and show that, after augmenting the covariance matrix by sufficient statistics dictated by a triangulation, its inverse is block graph-structured. BEGIN differs in that it identifies conditional independence exactly through the covariance conditions in Theorem~\ref{thm:begin}, without assuming strict positivity, a clique factorization, or a fixed parametric likelihood. In this sense, BEGIN can be viewed as a local building block for Markov random fields over binary or multinomial variables. Section~\ref{subsec:examples} provides examples illustrating how such local BEGIN structures can be assembled into larger Markov graphs.

We also note that when $\Supp(\bB) \subsetneq \{\pm 1\}^s$, the coefficients in part~(b) and the matrices in part~(c) need not be unique. The equivalence itself is unaffected: BEGIN is a structural statement about existence, factorization, and sparsity, not about unique representations.

\subsection{Examples}\label{subsec:examples}
\begin{figure}
\centering
\subfloat[BEGIN for $A \independent C \mid B$ as $\langle A,B \rangle \cap \langle B,C \rangle$.]{\includegraphics[width=.45\textwidth]{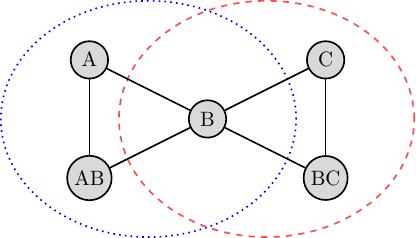}}
\hfill
\subfloat[BEGIN for a Markov chain $(A_1,A_2,A_3,A_4)$.]{\includegraphics[width=.5\textwidth]{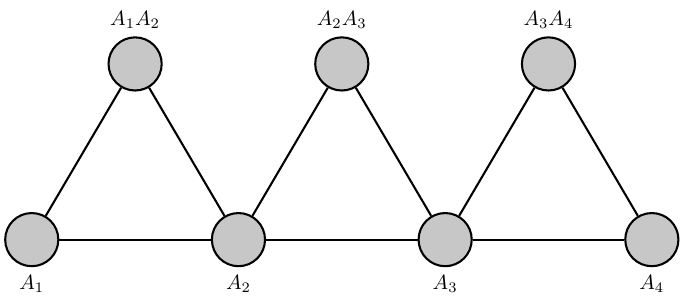}}\\[0.7em]
\subfloat[Left wing, center, and right wing of BEGIN for $B \independent (A_1,A_2,A_3) \mid (A_1A_2,A_2A_3,A_1A_3)$, corresponding to $\langle A_1A_2,A_1A_3,B \rangle \cap \langle A_1,A_2,A_3 \rangle$.]{\includegraphics[width=.80\textwidth]{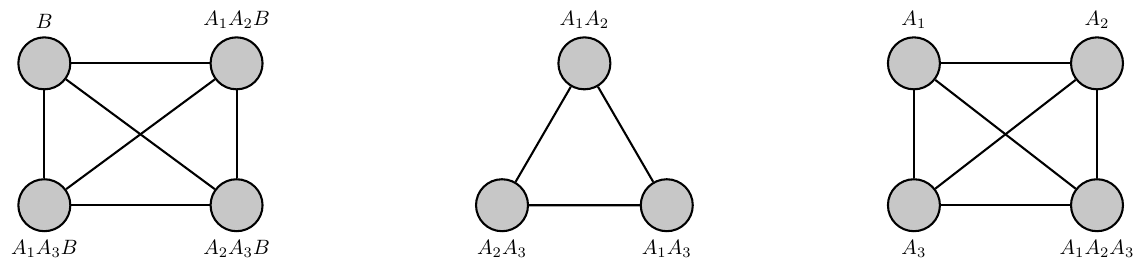}}
\caption{Examples of BEGIN, where conditional independence is represented through intersections of multiplicative groups of binary interactions.}
\label{fig:begin_examples}
\end{figure}

By incorporating interactions, the BEGIN framework can represent conditional independence structures that are difficult to display faithfully in classical ways and can serve as building blocks for Markov structures over binary variables.
\begin{enumerate}
\item \textbf{Three binary variables.} In the simplest case $r=s=t=1$, Figure~\ref{fig:begin_examples}(a) shows BEGIN for $A \independent C \mid B$. In addition to the original variables, BEGIN introduces the interaction nodes $AB$ and $BC$. The graph splits naturally into a left wing $\{A,AB\}$, a center $\{B\}$, and a right wing $\{C,BC\}$. The left wing together with the center generates $\langle A,B \rangle=\{1,A,B,AB\}$, the center together with the right wing generates $\langle B,C \rangle=\{1,B,C,BC\}$, and their intersection is $\langle B \rangle \setminus \{1\}=\{B\}$. 

\item \textbf{A binary first-order Markov chain.} Let $(A_1,\ldots,A_k) \in \{\pm 1\}^k$ be a (not necessarily stationary) first-order Markov chain. BEGIN contains the chain nodes $A_1,\ldots,A_k$ together with the interaction nodes $A_jA_{j+1}$ for $j=1,\ldots,k-1$; see Figure~\ref{fig:begin_examples}(b) for $k=4$. An unrestricted joint pmf on $\{\pm 1\}^k$ has $2^k-1$ free parameters, whereas a nonstationary first-order Markov model has only $2k-1$. BEGIN makes that reduction visible at the interaction-node level and suggests a sparse matrix representation for the chain. Moreover, Figure~\ref{fig:begin_examples}(b) is the union of the overlapping BEGIN molecules $\langle A_1,A_2 \rangle \cap \langle A_2,A_3 \rangle$ and $\langle A_2,A_3 \rangle \cap \langle A_3,A_4 \rangle$; more generally, the chain is assembled from the BEGIN molecules on $\langle A_j,A_{j+1} \rangle \cap \langle A_{j+1},A_{j+2} \rangle$, $j=1,\ldots,k-2$.

\item \textbf{A higher-order conditioning set.} \citet{brown2025belief} provide an example in which
\[
B \independent (A_1,A_2,A_3) \mid (A_1A_2,A_2A_3,A_3A_1).
\]
This form of conditional independence is not naturally expressed by a standard graph on $A_1$, $A_2$, $A_3$, and $B$. BEGIN, by contrast, yields a direct undirected graph on interaction nodes corresponding to $\langle A_1A_2,A_1A_3,B \rangle \cap \langle A_1,A_2,A_3 \rangle$; see Figure~\ref{fig:begin_examples}(c). 

\item \textbf{A Markov random field beyond the Ising model.} The BEGIN structure in the previous example can also serve as a building block for a four-node global Markov random field. Under the relabeling $X_1=B$, $X_2=A_1A_2$, $X_3=A_1A_2A_3$, and $X_4=A_1A_3$, Figure~\ref{fig:begin_examples}(c) represents $X_1 \independent X_3 \mid (X_2,X_4)$ through the group intersection $\langle X_1,X_2,X_4 \rangle \cap \langle X_2,X_3,X_4 \rangle.$ If $(X_1,X_2,X_3,X_4)\in\{\pm1\}^4$ further satisfies the BEGIN molecule $X_2 \independent X_4 \mid (X_1,X_3),$ then these two statements are exactly the nontrivial separation relations of the four-cycle standard graph $X_1-X_2-X_3-X_4-X_1$. Hence the distribution satisfies the global Markov property with respect to this graph. If the joint pmf is strictly positive, then $(X_1,X_2,X_3,X_4)$ is an Ising model on the four-cycle. Without strict positivity, the same pair of BEGIN molecules still defines a global Markov random field, but not necessarily an Ising model, because zeros in the joint pmf are allowed. Thus, this pair of BEGIN molecules yields a class of four-cycle global Markov random fields that is strictly larger than the Ising family.
\end{enumerate}

\subsection{The Hadamard prism}
The proof of Theorem~\ref{thm:begin}, provided in the supplementary material, relies on a linear mapping from $\RR^{2^p}$ to $\RR^{2^p \times 2^p}$ that packages the covariance algebra into a matrix operator. This mapping is closely related to convolution on $(\ZZ_2)^p$ and is diagonalized by the Walsh--Hadamard transform. Related constructions also appear in the literature on group-circulant matrices and Boolean Fourier analysis \citep{terras1999fourier, o2014analysis}. Because this operator is central to the covariance characterization underlying BEGIN and may also be of independent interest for future research, we give it a dedicated name.

\begin{definition}
For $\by \in \RR^{2^p}$, define the \emph{Hadamard prism} of $\by$ by
\begin{equation}\label{eq:hp}
\eta_p(\by)
:=
\frac{1}{2^p}\,\Hb_p\,\diag(\Hb_p\by)\,\Hb_p.
\end{equation}
\end{definition}

Because $\Hb_p$ is orthogonal up to scale, the eigenvalues of $\eta_p(\by)$ are proportional to the coordinates of $\Hb_p\by$. For binary interaction vectors, those coordinates are directly linked to cell probabilities and interaction means \citep{zhang2019bet}. The Hadamard prism also satisfies a recursion that is useful for second-moment calculations: for $\by_1,\by_2 \in \RR^{2^d}$,
\begin{equation}\label{eq:hp_recursion}
\eta_{d+1}\!\left(
\begin{pmatrix}
\by_1\\
\by_2
\end{pmatrix}
\right)
=
\begin{pmatrix}
\eta_d(\by_1) & \eta_d(\by_2)\\
\eta_d(\by_2) & \eta_d(\by_1)
\end{pmatrix}.
\end{equation}
This recursive form also suggests that the Hadamard prism may be useful beyond BEGIN for studying structured covariance patterns of binary variables.

\section{Approximating Conditional Independence for General Variables}\label{sec:approx}
The binary and multinomial theory above can be used as a multiresolution approximation device for general random vectors. Following \citet{zhang2019bet}, \citet{zhang2021beauty} and \citet{brown2025belief}, we consider binary expansion as a way to encode real-valued variables through data bits. The classical expansion
\[
U=\sum_{k=1}^{\infty} \frac{A_k}{2^k},
\qquad
U \in [-1,1],\quad A_k \in \{\pm 1\},
\]
suggests approximating $U$ by its $d$-bit truncation $U_d=\sum_{k=1}^d A_k/2^k$.

After marginal standardization, let $(\bU,\bV,\bW) \in [-1,1]^r \times [-1,1]^s \times [-1,1]^t$. Define the dyadic quantizer
\[
Q_d(x):=-1+2^{-d}+2^{1-d}\big\lfloor 2^{d-1}(x+1)\big\rfloor
\in \bigl\{-1+2^{-d},-1+2^{-d}+2^{1-d},\ldots,1-2^{-d}\bigr\},
\]
and apply it componentwise to vectors. The next result shows that exact dyadic conditional independence at every resolution implies the population notion and yields an explicit approximation rate under H\"older-type continuity.

\begin{theorem}\label{thm:beci}
For $(\bU,\bV,\bW) \in [-1,1]^r \times [-1,1]^s \times [-1,1]^t$ and $d\ge 1$, define
\[
\cU_d:=\sigma\bigl(Q_d(\bU)\bigr),
\qquad
\cV_d:=\sigma\bigl(Q_d(\bV)\bigr),
\qquad
\cW_d:=\sigma\bigl(Q_d(\bW)\bigr).
\]
Then the following statements hold.
\begin{enumerate}[(a)]
\item If $\cU_d \independent \cW_d \mid \cV_d$ for every $d\ge 1$, then $\bU \independent \bW \mid \bV$.

\item Suppose $\bU \independent \bW \mid \bV$. Assume there exist $\alpha \in (0,1]$ and constants $L_{\bU},L_{\bW}<\infty$ such that for all $\bv,\bv' \in [-1,1]^s$,
\[
\TV\bigl(\cL(\bU \mid \bV=\bv),\, \cL(\bU \mid \bV=\bv')\bigr)
\le
L_{\bU}\,\|\bv-\bv'\|_2^{\alpha},
\]
and
\[
\TV\bigl(\cL(\bW \mid \bV=\bv),\, \cL(\bW \mid \bV=\bv')\bigr)
\le
L_{\bW}\,\|\bv-\bv'\|_2^{\alpha}.
\]
Define
\[
\Delta_d:=\sup_{S\in\cU_d,\,T\in\cW_d}
\Eb\Bigl[\bigl|\Pb(S\cap T\mid \cV_d)-\Pb(S\mid \cV_d)\,\Pb(T\mid \cV_d)\bigr|\Bigr].
\]
Then, for every $d\ge 1$,
\[
\Delta_d \le L_{\bU}L_{\bW}\, s^{\alpha}\,2^{2\alpha(1-d)-2}.
\]
In particular, $\Delta_d \to 0$ as $d \to \infty$.
\end{enumerate}
\end{theorem}

Part~(a) shows that exact dyadic conditional independence at every resolution implies the population statement. Part~(b) controls the converse direction quantitatively. Without continuity assumptions, discretization can create spurious conditional associations or Simpson's paradox; see, for example, \citet{gong2021judicious}. Under H\"older-type regularity, however, the dyadic approximation error decays at an explicit rate. This provides theoretical support for using BEGIN on the leading data bits of continuous or mixed-type variables as a principled approximation to conditional independence.

\section{Discussion}\label{sec:discussion}
This note establishes an exact covariance-based graphical characterization of conditional independence for arbitrary multivariate binary data in the binary-expansion interaction basis, including singular multinomial encodings. The characterization is distribution-free and is expressed in objects that are natural from the perspective of multiresolution binary expansion. For multinomial and discretized continuous variables, the same viewpoint provides a principled way to relate exact bit-level statements to approximation results for more general variables. 

These results suggest several directions for future work. One concerns \emph{structure learning}: how should one estimate the sparse BEGIN graph efficiently from finite samples when the interaction feature space is large? The Schur-complement characterization and the BELIEF representation suggest nodewise procedures, regularized inverse problems, and screening rules tailored to interaction groups, but developing their finite-sample properties lies beyond the scope of this note. A second concerns \emph{statistical theory}: high-dimensional consistency, robustness to approximate sparsity, and finite-sample guarantees remain open. A third concerns \emph{causal and scientific interpretation}: the bit-level perspective suggests resolution-dependent notions of adjustment, mediation, and an approximation of causality.

\section*{Acknowledgments}
Zhang’s research was partially supported by NSF grants DMS-2152289 and TI-2449855, as well as BSF grant 2024055. The initial formulation and proof of Theorem~\ref{thm:begin} were completed while Zhou was a junior student at the Princeton International School of Mathematics and Science (PRISMS).  Zhou and Zhang thank PRISMS for supporting research collaborations involving high school students.  The Hadamard prism was developed during Zhang's visit to Michael Baiocchi at Stanford University.  Zhang thanks Baiocchi and Stanford University for the hospitality.  The authors thank Michael Baiocchi, Emmanuel Cand\`es, Peng Ding, Fang Han, Jan Hannig, Daniel Kessler, Han Liu, Yufeng Liu, Xiao-Li Meng, Heyang Ni, Art Owen, Evan Schwartz, Chengchun Shi, Daniel Yekutieli, Wan Zhang, Yuhao Zhou, Hongtu Zhu, and Jos\'e Zubizarreta for helpful comments and discussions.

\bibliography{KZ_BIB_121625.bib}
\bibliographystyle{chicago}
\end{document}